\documentclass[12pt]{article}

\usepackage{amsmath}
\usepackage{amssymb}
\usepackage{amscd}
\usepackage[all]{xy}

\newtheorem{theorem}{Theorem}[section]
\newtheorem{definition}[theorem]{Definition}
\newtheorem{corollary}[theorem]{Corollary}
\newtheorem{lemma}[theorem]{Lemma}
\newtheorem{proposition}[theorem]{Proposition}
\newtheorem{remark}[theorem]{Remark}

\newenvironment{proof}{\begin{trivlist} \item[]{\bf Proof.}}
{\par\hfill $\square$\end{trivlist}}

\newcommand{\C}{\mathbb{C}}

\newcommand{\R}{\mathbb{R}}

\newcommand{\N}{\mathbb{N}}

\newcommand{\ddc}{{{\rm dd}^{\rm c}}}
\newcommand{\Lone}{{{\rm L}^1}}
\newcommand{\Ltwo}{{{\rm L}^2}}

\newcommand{\PSH}{{\rm PSH}}

\renewcommand{\d}{{\rm d}}
\newcommand{\supp}{{\rm supp}}

\newcommand{\norm}[1]{\left\Vert#1\right\Vert}

\begin{document}
\pagestyle{plain}
\title{Lyapunov exponents and bifurcation current for polynomial-like maps}
\author{Ngoc-Mai Pham}
\maketitle
\begin{abstract}
We study holomorphic families of polynomial-like maps depending on a
parameter $s$. We prove that the partial sums of largest Lyapunov
exponents are plurisubharmonic functions of $s$. We also study their
continuity and introduce the bifurcation locus as the support of
bifurcation currents.
\end{abstract}

 \section{Introduction}
In this paper, we study the dependence of Lyapunov exponents on
parameters for polynomial-like maps in any dimension.

Recall that a \textit{polynomial-like map} is a proper holomorphic
map $f:U\to V$ where $U,V$ are open subsets of $\C^k$ and $U\Subset
V$. In particular, $f$ defines a ramified covering over $V$. The
degree $d_t$ of the covering is the \textit{topological degree} of
$f$, it is equal to $\sharp f^{-1}(z), z\in V$, counting
multiplicities. The family of polynomial-like maps is very large.
One checks easily that small perturbations of $f$ define also
polynomial-like maps.

In [DS1], Dinh-Sibony constructed for such a map the
\textit{equilibrium measure} $\mu$ as follows: if $\Omega$ is an
arbitrary smooth probability measure on $V$ then $\mu$ is the limit
of $d_t^{-n}(f^n)^*\Omega$ in the sense of measures. The measure
$\mu$ does not depend on the choice of $\Omega$. It is totally
invariant: $d_t^{-1}f^*\mu = f_*\mu =\mu$, of maximal entropy and is
mixing. Following Oseledec [Os], $f$ admits $k$ Lyapunov exponents
with respect to $\mu$ that we denote by
$$\chi_{1}\geq \chi_{2}\geq \cdots \geq \chi_{k}.$$
They satisfy the following inequality (see [DS1])
$$\chi_1 + \chi_2 +\cdots + \chi_k \geq \frac {1}{2}\log d_t.$$

The support $J$ of $\mu$ is called the \textit{Julia set} (of
maximal order) of $f$. It is contained in the boundary of the
\textit{filled Julia set}
$$
K:=\{z\in U: \,f^n(z)\in U \quad \text{for every } n\geq 0\}
=\bigcap_{n\geq 1}f^{-n}(V). $$

\noindent We recall that the measure $\mu$ \textit{maximizes the
plurisubharmonic (p.s.h. for short) moments} [DS1, Prop.3.2.6]. That
is, if $\nu$ is a totally invariant probability measure and
$\varphi$ is a p.s.h. function on $U$ then $\int \varphi \d \nu
\leq \int \varphi \d\mu.$

 Note that the study of
holomorphic endomorphisms of $\C P^k$ of degree algebraic $d\geq 2$
can be reduced to the study of some polynomial-like maps. Indeed, we
can lift these maps to $\C^{k+1}$ and the restrictions of the lifted
maps to a large ball $V\subset \C^{k+1}$ are polynomial-like maps.

Now, consider a holomorphic family of polynomial-like maps:
$$f_s:U_s \to V_s, s \in \Lambda.$$
More precisely, we have a holomorphic map $F:\mathcal{U}\to
\mathcal{V}$, where $\mathcal{U}\subset \mathcal{V}$ are open sets
in $\Lambda \times \C^k$ and $\Lambda$ is a connected complex
manifold of dimension $m$. Define
$$U_s=\mathcal{U}\cap (\{s\}\times
\C^k),\quad V_s=\mathcal{V}\cap (\{s\}\times \C^k).$$ We assume that
$U_s\Subset V_s$ and that the restriction of $F$ to $U_s$ defines a
polynomial-like map $f_s:U_s \to V_s.$ We often identify $U_s$ and
$V_s$ to open sets in $\C^k$. Observe that the topological degree
$d_t$ of $f_s$ does not depend on $s.$

 Let us denote the equilibrium
measure of $f_s$ by $\mu_s$, the Julia set by $J_s$ and the filled
Julia set by $K_s$. We order the Lyapunov exponents by
$$\chi_1(s)\geq \chi_2(s)\geq \cdots \geq \chi_k(s).$$

 \noindent In this paper, we prove that
$$L_p(s):=\chi_1(s)+\chi_2(s)+\cdots +\chi_p(s)$$
is p.s.h. on $s$ for $1\leq p\leq k$. The case where $p=k$ was
proved in [DS1].

  We define the \textit{bifurcation current} associated to
$(f_s)_{s\in \Lambda}$ by $$B_F:=\ddc L_k .$$
 This is a positive closed $(1,1)$-current on $\Lambda$. The support
 of $B_F$ is called the \textit{bifurcation
locus} of $(f_s)_{s\in \Lambda}$. Since $L_k$ is locally bounded
$(L_k\geq\frac{1}{2}\log d_t)$, we can consider the \textit{higher
degree bifurcation currents}: $$B_F^i:=B_F\wedge\ldots\wedge
B_F\quad (i \,\, \text{times})$$ and the \textit{higher degree
bifurcation locus} $\supp(B_F^i)$ for $0\leq i\leq m$ (see
Definition \ref{new11}). We also introduce in Section \ref{new6} the
\textit{total bifurcation current} $\hat{B}_F$ on $\mathcal{U}$ such
that $\pi_*(\hat{B}_F)=B_F$ where $\pi:\mathcal{U}\to \Lambda$ is
the canonical projection (see Definition \ref{new12}).

When the measure $\mu_{s_0}$ of $f_{s_0}$ is \textit{PLB} (i.e the
p.s.h. functions on $V_{s_0}$ are $\mu_{s_0}$-integrable) we show
that $L_k$ is continuous in a neighborhood of $s_0$. Note that the
property "$\mu_s$ PLB" is equivalent to the fact that some natural
dynamical degree $d_s^*$ of $f_s$ satisfies $d_s^*<d_t$.

We then study the stability of $(f_s)_{s\in \Lambda}$ in Section
\ref{new4}. In the case of dimension $1$, we are able to extend the
results of Ma\~n\'{e}-Sad-Sullivan (see [MSS]) to the case of
polynomial-like maps. In particular, we will prove that the
stability of $J_s$, i.e the continuity of $J_s$ in the Hausdorff
sense, implies that the unique Lyapunov exponent $\chi(s)$ defines a
pluriharmonic function on $\Lambda$ (see also [DM] for the case of
rational maps on $\C P^1$). In the case of higher dimension, we show
that if the critical set $\mathcal{C}_s$ of $f_s$ does not intersect
$J_s$ and if $\mu_s $ is PLB for $s\in \Lambda$ then $L_k$ is
pluriharmonic and $(f_s)_{s\in \Lambda}$ is stable. In this case,
the bifurcation locus is empty. The condition on $\mathcal{C}_s$ is
often easy to check.

 Note that in a recent work [BB], Bassanelli-Berteloot gave another
sufficient condition for $L_k$ to be pluriharmonic, for holomorphic
maps in $\C P^k$ (see also Remark \ref{new9}). A similar problem for
H\'{e}non maps was studied by Bedford-Lyubich-Smillie (see [BS],
[BLS]).

Observe that if $f:U\to V$ is a polynomial-like map then
$f:f^{-1}(V')\to V'$ is also polynomial-like for $U\subset V'\subset
V$. The problems on families of maps that we consider are
semi-local. Then we can assume to simplify the notation that $V_s=V$
for every $s$ and $\mathcal{V}=\Lambda \times V$ in Section
\ref{new6} and Section \ref{new5}.

The readers who are not familiar with the horizontal currents and
slicing theory, may consult our Appendix A or [DS2].

\section{Partial sums of Lyapunov exponents}\label{new6}

In this section, we study the partial sums of largest Lyapunov
exponents of polynomial-like maps. We first prove the following
useful result.

\begin{proposition}\label{new8}
Let $(f_s)_{s\in \Lambda}$ be a holomorphic family of
polynomial-like maps as above. Then there exists a horizontal
positive closed current $\mathcal{R}$ on $\Lambda \times V$ of
bidimension $(m,m)$ such that for every $s\in \Lambda$ the slice
$\langle \mathcal{R},\pi,s\rangle$ is equal to $\mu_s$ where
$\pi:\Lambda\times V\to \Lambda$ is the canonical projection
\end{proposition}

\begin{proof} Let $S$ be an arbitrary horizontal current of bidimension $(m,m)$ on $
\Lambda \times V$, define $S_s:=\langle S, \pi,s\rangle$ the slice
of $S$. Let $\vartheta$ be a smooth probability measure with compact
support in $V$. Consider the horizontal positive closed current
$S:=\pi_V^*(\vartheta)$ of bidimension $(m,m)$ on $ \Lambda \times
V$, where $\pi_V$ is the canonical projection of $\Lambda \times V$
on $V$. Define
$$S_n:=\frac{1}{d_t^n}(F^n)^*S.$$
We identify $S_s$ with $\vartheta$. Then
$S_{n,s}=\frac{1}{d_t^n}(f_s^n)^*\vartheta$ converges weakly to
$\mu_s$ for $s\in \Lambda$ (see [DS1]).

Since the masses of $S_n$ are locally uniformly bounded there exists
a subsequence $(i_n)_{n\in \N}$ such that $S_{i_n}$ converge to a
horizontal positive closed current $\mathcal{R}$ on $\Lambda\times
V$. By definition $\supp(\mathcal{R})\subset \bigcup_{s\in
\Lambda}K_s$, where $K_s$ is the filled Julia set of $f_s$. Let
$s_0$ be a fixed point in $\Lambda$. Let $U_0$ be a subset of
$U_{s_0}$ and $\Lambda_0$ be a small neighborhood of $s_0$ such that
$\bigcup_{s\in \Lambda_0} K_s\Subset U_0\Subset \bigcap_{s\in
\Lambda_0} U_s$. Consider a smooth p.s.h. function $\psi$ on a
neighborhood $\Lambda_0\times V_0$ and a continuous form $\Omega$ of
maximal degree with compact support in $\Lambda_0$. By formula
(\ref{4}) in the Appendix A, for every $n$, we have
$$\int_{\Lambda}S_{i_n,s} (\psi)\Omega(s)=
\langle S_{i_n}\wedge\pi^*(\Omega), \psi\rangle.$$ Then
\begin{eqnarray}\label{5}
\int_{\Lambda}\mu_s (\psi)\Omega(s)= \langle
\mathcal{R}\wedge\pi^*(\Omega),
\psi\rangle=\int_{\Lambda}\mathcal{R}_s (\psi)\Omega(s).
\end{eqnarray}
 Define $u(s):=\mathcal{R}_s(\psi)$ for $s\in \Lambda_0$.
  By Proposition \ref{r2.5}, $u$ is p.s.h. on $\Lambda_0$.

 We want to prove that $\mu_s(\psi)$ is also a p.s.h. function.
 Consider a sequence $(s_m)$ converging to $s_0$ such
 that $\mu_{s_m}\rightharpoonup \nu $. Hence $f^*_{s_0}\nu=d_t\nu$.
 Since $\psi$ is uniformly continuous on $\Lambda_0\times V_0$ then
 for $\epsilon > 0$ and $s$ close to $s_0$:
 $$\int{\psi(s,z)\d\mu_s}\leq \int{(\psi(s_0,z)+\epsilon) \d\mu_s}.$$
 Then
 \begin{align*}
 \limsup_{m \to \infty}\int{\psi(s,z)} \d\mu_{s_m} &\leq
 \limsup_{m\to
 \infty}\int{(\psi(s_0,z)+\epsilon)} \d\mu_{s_m}\\
&=\int{(\psi(s_0,z)+\epsilon)}\d\nu.
\end{align*}
Since the equilibrium measure $\mu_{s_0}$ maximizes p.s.h. moments
and $\epsilon$ is arbitrarily small, we get:
$$\limsup_{s\to s_0} \mu_s(\psi) \leq \int{\psi(s_0,z)}\d\nu \leq
\int{\psi(s_0,z)}\d\mu_{ s_0}=\mu_{s_0}(\psi).$$ Therefore,
$\mu_{s}(\psi)$ is upper semi-continuous.

 Define
\begin{align*}
\mu^N_s(\psi) :&= \int {\frac {1}{d_t^N}\psi(f^N_s)^*\vartheta}\\
&=\int {\frac{1}{d_t^N}(f^N_s)_*(\psi)\vartheta}.
\end{align*}
Recall here that $(f^N_s)_*(\psi)(z):=\sum_{f_s^N(w)=z}\psi(w)$
where the roots of the equation $f_s^N(w)=z$ are counted with
multiplicities. Since $\psi (s,z)$ is p.s.h., $(f^N_s)_*(\psi)(s,z)$
is a p.s.h. function of $(s,z)$. Hence $\mu^N_s(\psi)$ is p.s.h. on
$\Lambda_0$. Since $\mu_s$ is the limit of $(f^N_s)^*\vartheta$, as
$N \to \infty$, then $\mu^N_s(\psi)$ converges to $\mu_s(\psi)$, as
$N \to \infty$. Consequently, $\mu_s(\psi)$ is a p.s.h. function.

 The
equality (\ref{5}) is valid for all continuous form $\Omega$ of
maximal degree with compact support in $\Lambda_0$ then
$\mu_s(\psi)=u(s)$ a.e on $\Lambda_0$. But these functions are
p.s.h. hence $\mu_s(\psi)=u(s)=\mathcal{R}_s(\psi)$ for all $s\in
\Lambda_0$. We deduce that $\mu_s(\psi)=\mathcal{R}_s(\psi)$ for
$\psi$ smooth with compact support in $\Lambda_0\times V_0$. Indeed,
we can write $\psi=\psi_1 -\psi_2$ where $\psi_1$ and $\psi_2$ are
p.s.h. on $\Lambda_0\times V_0$. Therefore, $\mathcal{R}_s=\mu_s$
for every $s\in \Lambda$.
\end{proof}

The following theorem generalizes a result of [DS1].

\begin{theorem}\label{t2.1} Let $(f_s)_{s\in \Lambda}$ be a
holomorphic family of polynomial-like maps with topological degree
$d_t\geq 2$ as above. Then the function $L_p(s)=\chi_1(s)+ \cdots +
\chi_p(s)$ is p.s.h. on $\Lambda$ for all $1\leq p\leq k$. In
particular, $L_p$ is upper semi-continuous.
\end{theorem}

 Let $f:U\to V$ be a
polynomial-like map, $U\Subset V \Subset \C^k$. Let $\mu$ denote its
equilibrium measure, $J$ its Julia set and $\chi_i$ the Lyapunov
exponents with $\chi_1\geq\cdots\geq \chi_k$. Let $Df^n(z)$ denote
the differential of $f^n$ at $z$. This is a linear map from the
tangent space $T_z$ of $\C^k$ at $z$ to the one at $f^n(z)$.

 Consider an orthonormal family of vectors $\{e_1(z),
\ldots , e_p(z)\}$ in $T_z$ and the linear space $e^p(z)$ generated
by $e_1(z), \ldots , e_p(z)$ for $z \in f^{-n}(U)$. The volume
$\lambda (e^p(z), Df^n)$ of the parallelotope determined by the
vectors
$$Df^n(e_1(z)), \ldots ,
Df^n(e_p(z))$$ is called the \textit{coefficient of expansion} in
the direction $e^p(z)$. It depends only on $e^p(z)$, not on the
choice generators. If $e^{p+q}(z)=e^p(z)\oplus e^q(z)$, we have
$$\lambda(e^{p+q}(z), Df^n) \leq \lambda (e^p(z), Df^n)\,
\lambda(e^q(z), Df^n).$$ \quad A $p$-vector $v$ in the exterior
product space $\bigwedge^p T_z$ is \textit{simple} if it can be
written as $v=e_1\wedge\ldots\wedge e_p$ where $e_1,\ldots,e_p$ are
vectors in $T_z$. Simple $p$-vectors generate $\bigwedge^p T_z$. The
map $f$ induces a linear map

 $$\wedge^p Df^n(z): \bigwedge^pT_z \to
 \bigwedge^p T_{f^n(z)}$$
 which is defined by
 $$ \wedge^p Df^n(z)(v):= Df^n(e_{1}) \wedge
 \ldots \wedge Df^n(e_{p}).$$

 \noindent
 For $z \in J$ define the Lyapunov $p$-dimensional characteristic
 number of
$e^p(z)$ by:
$$\chi(e^p(z)) := \limsup_{n\to \infty}
 \frac{1}{n} \log \lambda (e^p(z), Df^n).$$
 \noindent
By [Os], there exists a regular subset $E$ of $J$ satisfying
$f(E)\subset E$, $\mu(E)=1$ such that for all $z\in E$, $1\leq p\leq
k$

$$\chi(e^p(z)) = \lim_{n\to \infty}
\frac{1}{n} \log \lambda (e^p(z), Df^n)$$
$$\chi_1+\cdots+\chi_p=\sup_{ e^p(z)}\chi(e^p(z)).$$
Define
 $\norm{\wedge^pDf^n(z)} := \sup_{ e^p(z)}
\lambda (e^p(z), Df^n)$. We deduce from the previous discussion that
for $1\leq p\leq k$ and $z\in E$,
\begin{align}\label{new13}
\chi_1+\cdots+\chi_p= \lim_{n\to \infty} \frac{1}{n}
\log\norm{\wedge^pDf^n(z)}
\end{align}
and
\begin{align}\label{new14}
 \norm{\wedge^pDf^{m+n}(z)} \leq
\norm{\wedge^pDf^n(z)} \norm{\wedge^pDf^m(f^n(z))}.
\end{align}

\bigskip
\noindent
 \textbf{Proof of Theorem \ref{t2.1}.}
Define $
\psi_{p,n}(s,z):=\frac{1}{n}\log\norm{\bigwedge^pDf^n_s(z)}$ and
$\varphi_{p,n}(s):=\int{\psi_{p,n}(s,z)\d\mu_s}$ for all $z\in
f_s^{-n}U_s$. It is clear that $ \psi_{p,n}(s,z)$ is p.s.h. on a
sufficiently small neighborhood of $\bigcup_{s\in \Lambda} K_{s}$.
By formula (\ref{new13}), we have
 $$L_p(s) = \lim_{n\to \infty}\varphi_{p,n}(s).$$

Proposition \ref{new8} and \ref{r2.5} imply that $\varphi_{p,n}$ is
p.s.h. on $\Lambda$. The inequality (\ref{new14}) implies that
$\{\varphi_{p,2^n}\}$ is a decreasing sequence of p.s.h. functions.
Then the limit $L_p(s)$ is also p.s.h. on $\Lambda$.
$\hfill{\square}$\bigskip

\begin{remark}\label{r2.3}\rm
Let $\mathcal{I}$ be a subset of $\{1,2,\ldots,k\}$ and
$L_{\mathcal{I}}(s):=\sum_{i\in \mathcal{I}}\chi_i(s)$. If
$L_{\mathcal{I}}$ is a p.s.h. function of $s$ for any holomorphic
family $(f_s)$ of polynomial-like maps then $L_{\mathcal{I}}$ is one
of the previous sums $L_1,L_2,\ldots,L_k$. To prove this, we can use
the following family of maps $f_s:\C^k \to \C^k$ by
\begin{align*}
 (z_1,z_2,\ldots,z_k)\mapsto
(f_{1,s}(z_1),f_{2,s}(z_2),\ldots,f_{k,s}(z_k)),
\end{align*}
where $(f_{i,s})$ is a holomorphic family of polynomial maps in
dimension 1 for $1\leq i\leq k$.
\end{remark}

\begin{definition}\label{new11}\rm Theorem \ref{t2.1} allows us to
define the following positive closed current associated to
$(f_s)_{s\in \Lambda}$
$$B_F:=\ddc L_k.$$ We call it the \textit{bifurcation current} and its
support the \textit{bifurcation locus}. Since $L_k\geq
\frac{1}{2}d_t$, we can define the \textit{higher degree bifurcation
current}
$$B_F^i:=B_F\wedge \ldots \wedge B_F\quad (i \quad \text{times})$$
and the \textit{higher degree
bifurcation locus} as support of $B_F^i$ for $1\leq i \leq m$.
\end{definition}

\begin{remark}\label{new10}\rm
In [DM], DeMarco considered holomorphic families of rational maps on
$\C P^1$. She proved that the bifurcation locus and the complement
of the \textit{stable set} coincide. The stable set is the largest
open subset of $\Lambda$ where the Julia set depends continuously on
the parameter.
\end{remark}

 We also have the following stronger "variation" of Theorem
\ref{t2.1}.

\begin{corollary}\label{c2.4}
Let $(f_s)_{s\in \Lambda}$ be a holomorphic family of
polynomial-like maps as above. Then
$$L_{p,\lambda}(s):=\max(\chi_1(s), \lambda)+ \max(\chi_2(s),\lambda)
+\cdots + \max(\chi_p(s), \lambda)$$ is a p.s.h. function on
$\Lambda$ for all $\lambda \in \R\cup \{-\infty \}$.
\end{corollary}
\begin{proof}
It is clear that
$$L_{p,\lambda}(s)=\max \{ L_p(s),
L_{p-1}(s)+\lambda,\ldots, L_1(s)+(p-1)\lambda, p\lambda\}.$$ The
corollary follows.
\end{proof}

We can chose $S$ in Proposition \ref{new8} such that $\mathcal{R}$
has support in $\partial\mathcal{K}$, where
$\mathcal{K}:=\bigcup_{s\in \Lambda} K_s$. Note that $K_s$ depends
upper semi-continuously on $s$ and $\mathcal{K}$ is closed. Using
Ces\`{a}ro means, we can construct a positive closed current $R$
such that $F^*(R)=d_tR$ and $R_s=\mu_s$. More precisely, $R$ is a
limit of a subsequence of $(\mathcal{S}_n)$, where
$\mathcal{S}_n:=\frac{1}{n}\sum_{i=1}^nS_n$. A horizontal positive
closed current $\mathcal{R}$ such that $\supp(\mathcal{R}) \subset
\partial\mathcal{K}$, $F^*(\mathcal{R})=d_t\mathcal{R}$ and
$\mathcal{R}_s=\mu_s$ is called \textit{equilibrium current}
associated to the family $(f_s)_{s\in\Lambda}$.

\begin{remark}\label{7}\rm
In the case $\dim V=1$, the equilibrium current $\mathcal{R}$ is
unique. More precisely, the function $u(s,w)=\int_{\{s\}\times
V}\log|w-z|\d \mu_s(z) $ is a p.s.h. potential for $\mathcal{R}$.
Let $S$ be a horizontal positive closed current with slice mass unit
such that $\supp (S_s)$ is not a polar set in $V$ for all $s\in
\Lambda$ then $\frac{1}{d_t^n}(F^n)^*S$ converges weakly to
$\mathcal{R}$.

\end{remark}

Define $J(s,z):=\det\mathrm{Jac}f_s(z)$ then $\ddc
(\log|J|)=[\mathcal{C}_F]$, where $\mathcal{C}_F$ is the critical
set of $F$. Let $\mathcal{R}$ be an equilibrium current of $(f_s)$
on $\Lambda \times V$. Note that
$\mathcal{R}_s(\log|J|)=\mu_s(\log|J|)=L_k\geq \frac{1}{2}d_t$.
Theorem \ref{new2} in the Appendix A implies that the current
$\log|J|\mathcal{R}$ is well defined.

\begin{definition}\label{new12}\rm
We define the \textit{total bifurcation current} associated to
$(f_s)_{s\in \Lambda}$ as
$$\hat{B}_F:=\ddc (\log|J|\mathcal{R})=[\mathcal{C}_F]\wedge
\mathcal{R}$$ and \textit{total bifurcation locus} as its support.
Then the current $\hat{B}_F$ is positive closed. Observe that
$B_F=\pi_*(\hat{B}_F)$. As an immediate consequence, if
$\partial{\mathcal{K}}\cap [\mathcal{C}_F]=\varnothing$ then the
bifurcation locus, total bifurcation locus are empty and $L_k$ is
p.s.h. (see also Theorem \ref{t4.3}).
\end{definition}

\begin{remark}\label{new9}\rm
If $F=(f_s)_{s\in \Lambda}$ is a holomorphic family of holomorphic
endomorphisms in $\C P^k$, we can find a horizontal $(1,1)$-current
$\tau$ on $\Lambda \times \C P^k$ such that the slices
$\langle\tau,\pi,s\rangle$ are the Green currents associated to
$f_s$. Moreover $\tau$ has local continuous potentials. Then we can
define \textit{total j-bifurcation currents} by
$$\hat{B}_F^{(j)}:=[\mathcal{C}_F]\wedge \tau^j\quad \text{and}\quad
B_F^{(j)}:=\pi_*(\hat{B}_F^{(j)}), \quad 1\leq j\leq k.$$ The
currents $\hat{B}_F^{(k)}$ and $B_F^{(k)}$ correspond to the
currents $\hat{B}_F$ and $B_F$ defined previously. These definitions
are the starting point of our study which was developed
independently by Bassanelli-Berteloot (see [BB]). Note that
Bassanelli-Berteloot also obtained a nice formula for the sum of all
the Lyapunov exponents of endomorphisms of $\C P^k$ which
generalizes formula obtained by DeMarco (see [DM]) in the one
variable case.
\end{remark}

\section{Continuity of the sum of the exponents}\label{new5}
In this section we give a necessary condition so that the sum of all
the Lyapunov exponents depends continuously on the parameter. We use
here a notion of measure PLB introduced by Dinh-Sibony in [DS1].

Recall that a positive measure $\nu$ with compact support in an open
set $U\subset \C^k$ is said to be \textit{PLB} in $U$ if p.s.h.
functions on $U$ are $\nu$-integrable.
 In other words, if $\varphi$ is p.s.h. on $U$ we have: $\int
 \varphi \d\nu > -\infty$.
 The following theorem is stronger than the fact that $\mu_s \rightharpoonup
 \mu_{s_0}$.
 \begin{theorem}\label{t3.1}
 Let $f_s$, $\mu_s$ be as above. Assume that $V$ is
 Stein and $\mu_{s_0}$ is PLB. Then $\mu_s(\varphi)\to \mu_{s_0}(\varphi)$,
  as $s\to s_0$ for all p.s.h. function
 $\varphi $ on $V$.
 \end{theorem}

 Note that in [DS1], the authors proved that "$\mu_s$ PLB" is stable
 under small perturbations. In the
case of dimension $1$, the equilibrium measure of any
polynomial-like map is PLB. In the context of Theorem \ref{t3.1},
$\mu_s$ is PLB if $s$ close enough to $s_0$.

Since the problem here is local for $s$, we can replace $V$ by a
small perturbation such that the boundary $\partial V$ is real
analytic. Then the boundaries of $U_{s}$ and $\mathcal{U}$ are also
real analytic for $ s\in \Lambda_0$. Choose a neighborhood
 $ \Lambda_0$ of $s_0$ in $\Lambda$ and a Stein open set $W$ with smooth
 boundary such that:
 $$\bigcup_{s\in \Lambda_0}f^{-2}_s(V)\Subset W \Subset \bigcap_{s\in
 \Lambda_0}U_s.$$

 The following proposition and lemmas 3.1- 3.4 are refinements of the results in
 [DS1]. We refer to that paper for the proofs.

\begin{proposition}\label{p3.2}Let p be a positive integer. Then

(i) For $s\in \Lambda_0$ the norm of the operator
$$\mathcal{L}_s:=d_t^{-1}(f_{s})_*: \PSH(W)\cap \Ltwo(W)\to
\PSH(U_s)\cap \Ltwo(U_s)$$
 is uniformly bounded by a constant $A$.

(ii) There exist a positive integer $n_0\geq 2$, a constant $0<c<1$
and a neighborhood $\Lambda_1$ of $s_0$, $\Lambda_1\subset
\Lambda_0$ such that if $\varphi $ is a p.s.h. function on $V$ then
$$\norm{\mathcal{L}_s^{nn_0}\ddc\varphi}_{U_{s_0}}\leq c^n\norm{\ddc\varphi}_
{U_{s_0}},$$
 for $s\in \Lambda_1$.
 \end{proposition}

 To prove Theorem \ref{t3.1}, we can replace
$f_s$ by $f_s^{n_0}$. Then there exists $\tilde{U}$ so that for
$s\in \Lambda_1$, we have
\begin{eqnarray}\label{3}
 \bigcup_{s\in \Lambda_1}U_s
 \Subset  \tilde{U} \Subset V;\,\,\norm{\mathcal{L}_s^{n}\ddc\varphi}_{\tilde{U}}
 \leq c^n\norm{\ddc\varphi}_
 {\tilde{U}}.
\end{eqnarray}

\begin{proposition}\label{p3.3}
Let $U\Subset V$ be two open sets in $\C^k$. Let
$(\nu_{\theta})_{\theta\in \Gamma}$ be a family of probability
measures supported in a compact set $K\subset U$. Suppose that there
is a constant $B>0$ such that for all p.s.h. function $\psi$ on U
and for all $\theta \in \Gamma$, we have:
$\norm{\psi}_{\Lone(\nu_{\theta})} \leq B\norm{\psi}_{\Ltwo(U)}$.
Then there exists a constant $0<b<1$ such that
$$\sup_U\psi \leq b\, \sup_V\psi$$
for all p.s.h. function
 $\psi$ on $V$ satisfying:
 $\int \psi \d\nu_{\theta}=0$ for at least one $\theta$.
\end{proposition}

  Let $H$ denote the subspace of pluriharmonic functions in $\Ltwo(W)$
  and $H^{\bot}$ the cone of p.s.h. functions orthogonal to $H$.
  Let $\varphi $
  be a p.s.h. function on $W$ and $\varphi =u +v$ with $u\in H$
  and $v\in H^{\bot}$ the canonical decomposition of $\varphi$. We also have
  $\mathcal{L}_s\varphi=\mathcal{L}_{1,s}u+\mathcal{L}_{2,s}v+\mathcal{L}_{3,s}v$, where
  $\mathcal{L}_{1,s}: H\to H$, $\mathcal{L}_{2,s}:H^{\bot}\to H$ and
  $\mathcal{L}_{3,s}:H^{\bot}\to H^{\bot}$ are canonical linear maps associated to $\mathcal{L}_s$.

   \noindent Following Proposition \ref{p3.2}, for $s \in \Lambda_0$, we have
   \begin{eqnarray}\label{2}
   \norm{\mathcal{L}_{2,s}},\norm{\mathcal{L}_{3,s}}\leq A\quad \text{for}\quad s\in \Lambda_0.
\end{eqnarray}
We have
   $$\mathcal{L}_s^n\varphi=\mathcal{L}^n_{1,s}u+ \mathcal{L}^{n-1}_{1,s}\mathcal{L}_{2,s}v+
   \mathcal{L}^{n-2}_{1,s}\mathcal{L}_{2,s}\mathcal{L}_{3,s}v+\cdots + \mathcal{L}_{2,s}
   \mathcal{L}_{3,s}^{n-1}v+ \mathcal{L}^n_{3,s}v.$$
The following lemma is a consequence of the solution of
$\partial\bar{\partial}$-equation in a Stein open set.
\begin{lemma}\label{l3.4}
There exists a constant $C>0$ such that for $s\in\Lambda_0$, $v\in
H^{\bot}$, we have
$$\norm{\mathcal{L}_{3,s}v}_{\Ltwo(W)}\leq C \norm{\ddc v}_{W}.$$
\end{lemma}

\noindent
 We have an easy following lemma.
 \begin{lemma}\label{l3.6}
 Let $K\Subset U\subset \C ^k$ then there exists a positive constant
 $A(K,U)$ such that
 for all $\varphi$ p.s.h. on U, we have
 $$\norm{\varphi}_{\Ltwo(U)}\geq A(K,U)\norm{\ddc\varphi}_K.$$
 \end{lemma}

 \begin{proof}
Let $\Phi$ be a positive form with compact support in $U$ which is
equal to $(\ddc \norm{z}^2)^{k-1}$ on a neighborhood of $K$. Then we
have $\norm{\ddc \varphi}_K\leq \langle\ddc
\varphi,\Phi\rangle=\int_U\varphi\ddc \Phi$. It is clear that there
exists a positive constant $A(K,U)$ such that $\int_U\varphi\ddc
\Phi \leq A(K,U)\norm{\varphi}_{\Ltwo (U)}$. Lemma \ref{l3.6}
follows.
\end{proof}
\begin{lemma}\label{l3.7}
Let $U_0$ be an open subset of V satisfying $\tilde{U}\Subset
U_0\Subset V$. Then there exists a constant $A_1>0$ such that for
$s\in \Lambda_1$ and $\varphi$ p.s.h. on $V$, we have
$$\norm{\varphi}_{\Lone(\mu_s)}\leq A_1\norm{\varphi}_{\Ltwo(U_0)}.$$
\end{lemma}

 \begin{proof}
Define $ b(s):=\int u \d\mu_s$, $b_{j}(s):=\int
\mathcal{L}_{2,s}\mathcal{L}_{3,s}^jv \d\mu_s$ and $$h_{n}(s):= b(s)
+b_{0}(s)+\cdots+b_{n-1}(s) .$$ It is
 prove that in [DS1] that the sequence $(h_{n}(s))$ converge to $\mu_s(\varphi)$.

 Inequality (\ref{3}) and Lemma \ref{l3.4} imply
\begin{eqnarray*}
 \norm {\mathcal{L}_{3,s}^jv}_{\Ltwo(W)}
 &\leq& C\norm {\ddc\mathcal{L}_{3,s}^{j-1}v}_{W}\\
 &=&C\norm {\ddc\mathcal{L}_s^{j-1}\varphi}_{W}\\
 &\leq& A_2 c^j\norm {\ddc\varphi}_{\tilde{U}},
\end{eqnarray*}
 where $A_2=C/c$ does not depend on $s $ and $\varphi$.

Since $\mathcal{L}_{2,s}\mathcal{L}_{3,s}^jv$ is pluriharmonic (for
each $s$ fixed), we deduce from the last inequality and the
inequality (\ref{2}) that
$$|b_{j}(s)|\leq A_3 c^j\norm{\ddc\varphi}_{\tilde{U}}\quad \text{and}\quad |b(s)|\leq A_3
\norm{\varphi}_{\Ltwo(\tilde{U})},
$$
where $A_3$ is independent of $s$ and $\varphi$. The second
inequality is a consequence of the pluriharmonicity of $u$. By Lemma
\ref{l3.6} and the inequalities above, there exists a constant
$A_4>0$ such that
\begin{eqnarray}\label{new15}
\mu_s(\varphi)\geq -A_4\norm{\varphi}_{\Ltwo(U_0)}.
\end{eqnarray}
Define $\varphi^+:=\max(\varphi,0)$ then
$$\int|\varphi|\d \mu_s=\int(-\varphi + 2\varphi^+)\d \mu_s\leq
\big| \mu_s(\varphi)\big|+ 2\sup_W\varphi^+.$$ The submean
inequality implies that $\mu_s(\varphi)\leq \sup_W\varphi^+\leq
A_5\norm{\varphi}_{\Ltwo(U_0)}$, where $A_5$ is independent of $s$
and $\varphi$. By inequality (\ref{new15}), we have
$$\norm{\varphi}_{\Lone(\mu_s)}\leq A_1\norm{\varphi}_{\Ltwo(U_0)},$$
where $A_1:=\max(A_4,3A_5)$ is independent of $s$ and $\varphi$.
\end{proof}

\noindent \textbf{Proof of Theorem \ref{t3.1}.}  Let $V_0$ be an
open subset of $V$ such that $U_0\Subset V_0\Subset V$. By Lemma
\ref{l3.7}, the family $(\mu_s)_{s\in \Lambda_1}$ satisfies the
hypothesis of Proposition \ref{p3.3}. Then there exists $ 0<c_0<1$,
(independent of $s\in \Lambda_1$ and $\varphi$), such that

 $$
\sup_{V_0}(\mathcal{L}_s^{n+1}\varphi-\mu_s(\varphi))\leq c_0
\sup_{V_0} (\mathcal{L}_s^n\varphi -\mu_s(\varphi)).$$ Hence
 $$ \sup_{V_0}\mathcal{L} _s^n\varphi-\mu_s(\varphi)\leq
c_0^{n-1}(\sup_{U_0}\varphi -\mu_s(\varphi)).$$ By Lemma \ref{l3.7},
there exists a constant $M>0$ (independent of $s\in \Lambda_1$ and
$\varphi$) such that
$$ \sup_{V_0}\mathcal{L}_s^n\varphi -\mu_s(\varphi)\leq M  \norm{\varphi}_{\Ltwo(V_0)}c_0^n.$$
Since $W\Subset V_0$ and
$\mu_s(\mathcal{L}_s^n\varphi)=\mu_s(\varphi)$, we obtain
\begin{eqnarray}\label{new16}
0\leq \sup_{W}\mathcal{L}_s^n\varphi -\mu_s(\varphi)\leq M
\norm{\varphi}_{\Ltwo(V_0)}c_0^n.
\end{eqnarray}
Then
$$|\mu_s(\varphi)-\mu_{s_0}(\varphi)|\leq
|\sup_{W}\mathcal{L}_s^n\varphi-\sup_W \mathcal{L}_{s_0}^n\varphi|+
M  \norm{\varphi}_{\Ltwo(V_0)}c_0^n.$$

 Define $\mathcal{M}_{n}(s):= \sup_W\mathcal{L}_s^n\varphi$.
 By the last inequality, if $\mathcal{M}_{n}(s)$ is a function continuous at $s_0$ for every $n
$, then $\mu_s(\varphi)\to \mu_{s_0}(\varphi)$. Fin an index $n$. We
will prove the continuity of $\mathcal{M}_{n}(s)$ at $s_0$.

Let $(s_m)\to s_0$ and $z_m\in W$ so that
$\mathcal{L}^n_{s_m}\varphi(z_m)\geq \mathcal{M}_{n}(s_m)-1/m$. By
extracting a subsequence, we can assume that $z_m\to z_0\in
\overline{W}$. Since $\mathcal{L}^n_s(\varphi)$ is a p.s.h. function
of $(s,z)$, by upper semi-continuity property, we have $\limsup
\mathcal{L}^n_{s_m}\varphi(z_m)\leq
\mathcal{L}^n_{s_0}\varphi(z_0)$. Observe that
$\sup_{\overline{W}}\psi=\sup_{W}\psi$ for all $\psi$ p.s.h. on $V$
since $W$ has smooth boundary. Then
$\mathcal{L}^n_{s_0}\varphi(z_0)\leq \mathcal{M}_{n}(s_0)$. Hence
$\limsup_{s\to s_0} \mathcal{M}_{n}(s)\leq \mathcal{M}_{n}(s_0)$.

 Fixed a positive number
$\epsilon$ and a no-critical value $x_0\in W$ of $f^n_{s_0}$ such
that $\mathcal{L}^n_{s_0}\varphi(x_0)\geq
\mathcal{M}_{n}(s_0)-\epsilon$. If $r>0$ is small enough then
$B(x_0,2r)$ is contained in $W$ and does not intersect the set of
critical values of $f^n_s$ for $s$ close enough to $s_0$. We see
that $\mathcal{L}^n_s\varphi$ converges to
$\mathcal{L}^n_{s_0}\varphi$ in $\Lone(B_r)$, where we denote
$B_r:=B(x_0,r)$.

\noindent Hence when $s\to s_{0}$, we have

$$\frac{1}{\text{vol}(B_r)}\int_{B_r} \mathcal{L}^n_s\varphi(x)\d x
\to \frac{1}{\text{vol}(B_r)}\int_{B_r}
\mathcal{L}^n_{s_0}\varphi(x)\d x .$$
 On the other hand, the submean inequality gives us
\begin{eqnarray*}
\sup_{B_r}\mathcal{L}^n_s(\varphi)&\geq&
\frac{1}{\text{vol}(B_r)}\int_{B_r} \mathcal{L}^n_s\varphi(x)\d x\\
&\geq&\frac{1}{\text{vol}(B_r)}\int_{B_r}
\mathcal{L}^n_{s_0}\varphi(x)\d x-\epsilon\\
 &\geq&
\mathcal{L}^n_{s_0}\varphi(x_0)-\epsilon \geq
\mathcal{M}_{n}(s_0)-2\epsilon,
\end{eqnarray*}
for $s$ close enough to $s_0$.
 Therefore
$$\liminf_{s\to s_0}\mathcal{M}_{n}(s)\geq \liminf_{s\to s_0}
\sup_{B_r}\mathcal{L}^n_s(\varphi)\geq
\mathcal{M}_{n}(s_0)-2\epsilon.$$ It follows that $\liminf_{s\to
s_0}\mathcal{M}_{n}(s)\geq \mathcal{M}_{n}(s_0)$. Therefore
$\mathcal{M}_{n}(s)$ is continuous at $s_0$.

$\hfill{\square}$\bigskip

 It is well-known that the
Lyapunov exponent is continuous in the space of rational functions
on $\C P^1$ (see [Ma]). We have the following result for families of
polynomial-like
 maps.
 \begin{theorem}\label{t3.9}
 Let $(f_s)_{s\in \Lambda}$ be a holomorphic family of
 polynomial-like
 maps as above. If $\mu_{s_0}$ is PLB  and $V$ is Stein then
 the sum $L_k(s)$ of all the Lyapunov exponents of $f_s$ is continuous on
 a neighborhood of $s_0$.
\end{theorem}
\begin{proof}\quad Because $\mu_s$ is PLB in a small neighborhood
of $s_0$ then it is sufficient to prove that $L_k$ is continuous
at $s_0$. Define $\varphi_s:=\log|\det \mathrm{Jac}(f_s)|$.
Replace $V$ by a Stein open subset of $U_{s_0}$ then we can assume
that $\varphi_s$ is p.s.h. on $V$. This function is continuous on
$(s,z)$ with value
 in $[-\infty,\infty[$. Since $\partial W$ is smooth,
 $\sup_{W}\mathcal{L}_s^n\varphi_s$ is continuous for every $n$.
 From inequality (\ref{new16}), we deduce that
  $L_k(s)=\mu_s(\varphi_s)$ is continuous at $s_0$.
  \end{proof}
In the case of dimension 1, every family of polynomial-like maps
satisfies the hypothesis of Theorem \ref{t3.9}. We have the
following corollary.
\begin{corollary}\label{new17}
Let $(f_s)_{s\in \Lambda}$ be a holomorphic family of
 polynomial-like maps as above in dimension one. Then
 the unique Lyapunov exponent $\chi(s)$ of $f_s$ is continuous.
\end{corollary}

  We also obtain the following corollary.

\begin{corollary}\label{c3.10}
Let $\{f_s:\C P^k\to\C P^k\}_{s\in \Lambda}$ be a holomorphic family
of holomorphic endomorphism of $\C P^k$ of algebraic degree $d\geq
2$. Then the sum $L_k$ of all the Lyapunov exponents of $f_s$ is
continuous.
\end{corollary}

\begin{proof} We can, locally on $\Lambda$, lift $(f_s)_{s\in \Lambda}$ to a holomorphic family
of homogeneous polynomials $\{F_s:\C^{k+1}\to \C^{k+1}\}_{s\in
\Lambda}$. Then by Theorem \ref{t3.9}, the sum $\tilde{L}_{k+1}(s)$
of all the Lyapunov exponents of $F_s$ is continuous on $\Lambda$.
Hence, $L_k(s) =\tilde{L}_{k+1}(s)-\log d$ is also continuous.
\end{proof}

\section{Stability of the Julia sets}\label{new4}

The purpose of this section is to find some sufficient conditions so
that the family of polynomial-like mappings $\{f_s: U_s\to
V_s\}_{s\in \Lambda}$ is stable and the sum of all Lyapunov
exponents of $f_s $ is a pluriharmonic function. We say that
$(f_s)_{s\in \Lambda}$ is stable if the Julia set $J_s$ depends
continuously on $\Lambda$ in the Hausdorff sense.

The stability of the Julia set for rational maps has been studied by
Ma\~n\'{e}-Sad-Sullivan [MSS] (see also [Mc], [DH], [DM]). Their
results can be extended to the case of polynomial-like maps. We have
the following result.

\begin{proposition}\label{p4.1}
Let $\{f_s: U_s\to V_s\}_{s\in \Lambda}$ be a holomorphic family of
polynomial-like maps in dimension one. Let $s_0$ be a point in
$\Lambda$. Then the following conditions are equivalent:

 (1) The number of
attracting cycles of $f_s$ is locally constant at $s_0$.

 (2) The
maximum period of attracting cycles of $f_s$ is locally bounded at
$s_0$.

(3) The Julia set moves holomorphically at $s_0$.

 (4) For $s$
sufficiently close to $s_0$, every periodic point of $f_s$ is
attracting, repelling or persistently indifferent.

 (5) The Julia set
$J_s$ depends continuously on $s$ (in the Hausdorff topology) on a
neighborhood of $s_0$.

 Suppose in addition that there are holomorphic maps $c_i:\Lambda
\to \C$ which parameterize the critical points of $f_s$. The the
following condition is also equivalent to those above:

 (6) There is a neighborhood $U$ of $s_0$ such that for $s$ in $U$,
 $c_i(s)\in J_s$ if and only if $c_i(s_0)\in J_{s_0}$.
\end{proposition}

 We also have the following theorem in the case of dimension 1.

\begin{theorem} \label{t4.2}Let $\{f_s: U_s\to V_s\}_{s\in \Lambda}$
be a holomorphic family of polynomial-like maps with topological
degree $d_t\geq 2$. If $(f_s)_{s\in \Lambda}$ stable in $\Lambda$
then the unique Lyapunov exponent $\chi(s)$ is a pluriharmonic
function.
\end{theorem}

\begin{proof}
Let $N(s)$ denote the number of critical points of $f_s$ counted
without multiplicity. Let $D(f)$ denote the set of $s'\in \Lambda$
such that $N(s)$ does not have a local maximum at $s'$. This is a
proper subvariety of $\Lambda$.

 If $s_0\notin D(f)$
then there is a neighborhood $\Lambda_0$ of $s_0$ in $\Lambda$ and
holomorphic functions $c_j:\Lambda_0\to \C $, $j=1,2,\ldots,N$,
parameterizing the critical points of $f_s$ (counted with
multiplicity).

 Therefore, $f_s'(z)=\Pi_{j=1}^{N}(z-c_j(s))h_s(z)$,
where $h_s(z)$ is a holomorphic function of $(s,z)$ which does not
vanish. Hence,
$$\chi(s)=\sum_{j=1}^{N}\int \log|z-c_j(s)|\d\mu_s + \int
{\log|h_s(z)|\d\mu_s}.$$ By Propositions \ref{new8} and \ref{r2.5},
$\int {\log|h_s(z)|\d\mu_s}$ is a pluriharmonic function of $s$. We
want to prove that for each $j$, the function
$$\lambda_j(s):=\int \log|z-c_j(s)|\d\mu_s $$
is pluriharmonic.

Fix an index $j$. Using a suitable holomorphic change of coordinate
$(s,z)$, we can assume that: $c_j(s)=0$ for $s\in \Lambda_0$. Then
$\lambda_j(s)=\int \log|z|\d\mu_s$. Let $A_{s,n}$ be the set of
repelling periodic points $a_{n,i}(s)$ of period $n$. It follows
from Proposition \ref{p4.1} that $I_n=\sharp A_{s,n}$ is independent
of $s$ and $a_{n,i}(s)$ is a holomorphic function of $s$, for $1\leq
i \leq I_n$. Define
$$\mu_{s,n}:= \frac{1}{d_t^n} \sum ^{I_n}_{i=1} \delta_{ a_{n,i}(s)}.$$
By Proposition \ref{p4.1}, either $0\in J_s$ for every $s\in
\Lambda_0$ or $0 \notin J_s$ for every $s\in \Lambda_0$.
  Let $\overline{\triangle(r)}$ denote the closed disk of center $0$ and of radius
$r$ with $r$ small. We define for each positive fixed number $r$,
\begin{eqnarray*}
I_{n,r}&:=&\big\{i: a_{n,i}(s_0)\notin \overline{\triangle(r)} \big\}\\
\mu_{s,n,r}&:=&\frac {1}{d_t^n}\sum_{i\in I_{n,r}}\delta
_{a_{n,i}(s)}\\
\mu_{s_0,r}&:=&\mu_{s_0}|_{J_{s_0}\backslash
\overline{\triangle(r)}}\\
\lambda_{j,r}(s)&:=&\limsup_{n\to \infty}\int \log |z|\d\mu
_{s,n,r}.
\end{eqnarray*}

By Proposition \ref{p4.1}, there exist holomorphic motions $\{\phi
_s: J_{s_0} \to \C\}_{s\in \Lambda_0}$ such that: $\phi_s \circ
f_{s_0}(z)=f_s \circ \phi_s(z)$, $\phi _s(J_{s_0})=J_s$, $\phi
_s(a_{n,i}(s_0))=a_{n,i}(s)$ and if $0\in J_s$ then $\phi_s(0)=0$.
We have
$$a_{n,i}(s)\in J_s \backslash \phi_s(J_{s_0}\cap \overline
{\triangle (r)})\quad \text {for all }i \in I_{n,r}.$$ It follows
that $\int \log|z|\d \mu_{s,n,r}$ is pluriharmonic on $s\in
\Lambda_0$ and bounded by a constant independent of $n$. Put
$\mu_{s,r}:=\mu_s|_{J_s \backslash \phi_s(J_{s_0}\cap \overline
{\triangle (r)})}$. We have $\mu_{s,n,r}\rightharpoonup \mu_{s,r}$
as $n\to \infty.$ These measures have support out of a neighborhood
of $0$. Hence $\int \log|z|\d \mu_{s,n,r}$ converges to $\int
\log|z|\d \mu_{s,r}$. This implies that $\lambda_{j,r}(s)$ is
pluriharmonic on $s\in \Lambda_0$. When $r>0$ is small and decreases
to $0$, $\lambda_{j,r}(s)$ decreases to $\lambda_j(s)$. Hence
$\lambda_j(s)$ is pluriharmonic on $\Lambda_0$. Therefore, $\chi(s)$
is a pluriharmonic function on $\Lambda \backslash D(f)$. By
Corollary \ref{new17}, $\chi(s)$ is continuous on $\Lambda $. It
implies that $\chi(s)$ is pluriharmonic on $\Lambda$.
\end{proof}

The following result is valid in any dimension.
\begin{theorem}\label{t4.3}
Let $\{f_s:U_s \to V_s \}_{s\in \Lambda }$ be a holomorphic family
of polynomial-like maps of topological degree $d_t\geq 2$ and
$\mathcal{C}_s$ the critical set of $f_s$. Assume that $\mu_s$ is
PLB and $\mathcal{C}_s \cap J_s=\varnothing$ for $s\in \Lambda$.
Then

 (i) The sum $L_k(s)$
of all the Lyapunov exponents of $f_s$ is a pluriharmonic function.
In particular, the bifurcation locus is empty.

 (ii) The family $(f_s)_{s\in \Lambda}$ is stable.
\end{theorem}

Note that a polynomial-like mapping satisfying the condition
$\mathcal{C}_s \cap J_s=\varnothing$ is, in general, not uniformly
hyperbolic on $J_s$.

 By [DS1], $f_s$ admits repelling periodic points on $J_s$. We have
 the following lemma, see [FS1] for the proof.

\begin{lemma}\label{l4.4}
Let $\{f_s:U_s \to V_s \}_{s\in \Lambda }$ be a holomorphic family
of polynomial-like maps of topological degree $d_t\geq 2$. Then for
all $s_0\in \Lambda $, there exists a neighborhood $\Lambda_{s_0}$
of $s_0$, a positive integer $N$ and repelling periodic points
$p(s)\in J_s$ such that $p(s)$ depends holomorphically on
$s\in\Lambda_{s_0}$ and $f^N_s(p(s))=p(s)$.
\end{lemma}

\noindent
 \textbf{Proof of Theorem \ref{t4.3}.}
 (i) Let $\mathcal{E}_s$ denote the exceptional set of $f_s$, i.e
 the set of point $z\in V_s$ such that the measure
 $d_t^{-n}\sum_{f_s^n(w)=z}\delta_w$
 does not converge to $\mu_s$. Since $\mu_s$ is PLB then $\mathcal{E}_s$
 is contained
 in the postcritical set $\bigcup_{n\geq 1} f_s^n(C_s)$ of $f_s$ for all
 $s\in \Lambda$ (see [DS1]). Hence $\mathcal{E}_s \cap
 J_s=\varnothing$. Fix a point $s_0$ in $\Lambda$ and let $p(s)$
 be as in Lemma \ref{l4.4}. Define
 $$\mu_{s,n}=\frac {1}{d_t^n}\sum _{i=1}^{d_t^n}\delta
 _{p_{n,i}(s)}$$
 where $p_{n,i}(s)$ are preimages of $p(s)$ by $f_s^n$.
Since $J_s$ is invariant by $f_s^{-1}$, all the points $p_{n,i}(s)$
are in $J_s$. The condition ${\mathcal{C}_s \cap J_s}=\varnothing$
implies that $p_{n,i}(s)$ depends holomorphically on $\Lambda_{s_0}$
and $\log |\det\mathrm{Jac}(f_s)|$ is pluriharmonic on a
neighborhood of $J_s$. These and the property that
$\mu_{s,n}\rightharpoonup \mu_s$ imply
\begin{align*}
L_k(s)&=\int \log |\det \mathrm{Jac}(f_s)|\d\mu_s \\
&=\lim_{n\to \infty}\int \log |\det \mathrm{Jac}(f_s)|\d\mu_{s,n}\\
&=\frac{1}{d_t^n}\lim_{n\to \infty}\sum_{i=1}^{d_t^n}\log |\det
\mathrm{Jac}(f_s)(p_{n,i}(s))|.
\end{align*}
 Therefore
$L_k(s)$ is pluriharmonic.\bigskip

 (ii) Observe that the family of holomorphic maps
$p_{n,i}:\Lambda_{s_0}\to W$ is normal where $W$ is an open set such
that $J_s\subset W\Subset \C^k$ for $s\in \Lambda_{s_0}$. Consider
the family $\mathcal{F}$ of all the maps $v:\Lambda_{s_0}\to W$ that
we obtain as limit, locally uniformly on $\Lambda_{s_0}$, of a
subsequence of $(p_{n,i})$. Hence $\mathcal{F}$ is a normal family
and $ \bigcup_{v\in \mathcal{F}}\upsilon (s)=J_s$ since
$\mu_{s,n}\rightharpoonup \mu_s$.
 It follows that
$(f_s)_{s\in \Lambda}$ is stable.$\hfill \square$

\begin{corollary}
Let $\{f_s:\C P^k \to \C P^k \}_{s\in \Lambda }$ be a holomorphic
family of holomorphic endomorphisms of algebraic degree $d\geq 2$.
Assume that $\mathcal{C}_s \cap J_s=\varnothing$ for $s\in \Lambda$,
where $\mathcal{C}_s$ is the critical set of $f_s$. Then

  (i) The sum $L_k(s)$
of all the Lyapunov exponents of $f_s$ is pluriharmonic on
$\Lambda$.

 (ii) The family $(f_s)_{s\in \Lambda}$ is stable.
\end{corollary}
\begin{proof} We can, locally on $\Lambda$, lift $(f_s)_{s\in \Lambda}$ to a holomorphic family
of homogeneous polynomials $\{F_s:\C^{k+1}\to \C^{k+1}\}_{s\in
\Lambda}$. Since $\mathcal{C}_{s}\cap J_s=\varnothing$, it implies
that $\tilde{\mathcal{C}}_{s}\cap \tilde{J}_{s}=\varnothing$ where
$\tilde{\mathcal{C}}_s$ and $\tilde{J}_s$ are the critical set and
the Julia set of $F_s$ respectively. By Theorem \ref{t4.3}, the sum
$\tilde{L}_{k+1}(s)$ of all the Lyapunov exponents of $F_s$ is
pluriharmonic on $\Lambda$ and $(F_s)_{s\in \Lambda}$ is stable.
Hence, $L_k(s) =\tilde{L}_{k+1}(s)-\log d$ is also pluriharmonic.

 If $\pi:\C^{k+1}\backslash \{0\} \to \C P^k$ denotes the
canonical projection, we have $J_{s}=\pi(\tilde{J}_{s})$. Note that
there is a neighborhood of $0$ in $\C ^{k+1}$ which does not
intersect $\tilde{J}_s$. Since $(F_s)_{s\in \Lambda}$ is stable then
$(f_s)_{s\in \Lambda}$ is also stable.
\end{proof}

\begin{appendix}
\section{Appendix: horizontal currents}

 Let $\Lambda$ and $V$ be two
bounded open subsets of $\C^m$ and $\C^k$ respectively. Let $\pi$
and $\pi_V$ denote the canonical projections of $\Lambda \times V$
on $\Lambda$ and on $V$. Let $R$ be a positive closed current of
bidegree $(k,k)$ on $\Lambda \times V$. We say that $R$ is {\it
horizontal} if $\pi_V(\supp(R))\Subset V$. Dinh-Sibony (see [DS2])
proved that the slice measure $\langle R,\pi,s\rangle$ is defined
for every $s$ and its mass is independent of $s$. We call this mass
the {\it slice mass} of $R$. We can consider $\langle
R,\pi,s\rangle$ as the intersection (wedge-product) of the current
$R$ with the current $[\pi^{-1}(s)]$ of integration on
$\pi^{-1}(s)$. The slice measure are characterized by the following
formula:

\begin{eqnarray}\label{4}
\int_{\Lambda}\langle R,\pi,s\rangle (\psi)\Omega(s)= \langle
R\wedge \pi^*(\Omega), \psi\rangle,
\end{eqnarray}
for every continuous $(m,m)$-form $\Omega$ with compact support in
$\Lambda$ and every continuous test function $\psi$ on
$\Lambda\times V$. Note that the formula (\ref{4}) is valid in the
general case where $\pi$ is an holomorphic submersion between two
complex manifolds $\mathcal{V}$ and $\Lambda$.

More over, we have the following proposition.

\begin{proposition}\label{r2.5}
Let $R$ be a horizontal positive closed current on $\Lambda \times
V$ and $\psi$ be a p.s.h. function on a neighborhood of $\supp(R)$.
Then the function
$$\varphi(s):=\int \psi(s,\cdot)\langle R,\pi,s\rangle$$
is p.s.h. on $\Lambda$ or equal to $-\infty$ identically.
\end{proposition}

\begin{proof}
This Proposition is a consequence of [DS2, Theorem 2.1], where the
authors consider the case of continuous p.s.h. function $\psi$. We
obtain the general case by using a sequence of smooth p.s.h.
functions decreasing to $\psi$.
\end{proof}

We refer to [De], [FS2] and [DS2] for the theory of intersection of
currents. We now prove the following theorem.

\begin{theorem}\label{new2}
Let $R$ be a horizontal positive closed current. Let $u$ be a p.s.h.
function on $\Lambda \times V$. Assume there exists $s_0\in \Lambda$
such that $\langle R, \pi ,s_0\rangle(u)\ne-\infty$. Then the
current $uR$ has locally finite mass in $\Lambda \times V$. In
particular, the positive closed current $\ddc u \wedge R$ is well
defined.
\end{theorem}

\begin{proof}
Consider open sets $\Lambda_0\Subset \Lambda$, $V_0\Subset V$ such
that $s_0\in \Lambda_0$ and $R$ is a horizontal positive closed
current on $\Lambda_0\times V_0$. Since the problem is local then it
is sufficient to prove that $uR$ has locally finite mass in
$\Lambda_0\times V_0$.

Let $A$ be a matrix of size $m\times k$ with complex coefficients.
We can consider $A$ as a point in $\C ^{mk}$. Define the affine map
$H:\C ^{mk}\times\C^m\times\C^k \to \C^m\times\C^k$ by
$$H(A,s,z):=(s-Az,z).$$
There exists a small ball $B(0,r)$ in $\C^{mk}$ such that
$\mathcal{R}:=H^*(R)$ is a horizontal positive closed current on
$\tilde{\Lambda}\times V_0$, where $\tilde{\Lambda}:=B(0,r)\times
\Lambda_0$.

 Define
$\tilde{u}:\tilde{\Lambda}\times V_0\to \R \cup \{-\infty\}$ by
$\tilde{u}(A,s,z):=u(s+Az,z)$. Then $\tilde{u}$ is p.s.h. on
$\tilde{\Lambda}\times V_0$. Consider $v:B(0,r)\times V_0\to \R \cup
\{-\infty\}$ by $v(A,s):=\langle
\mathcal{R},\pi_{\tilde{\Lambda}},(A,s)\rangle(\tilde{u})$, where
$\pi_{\tilde{\Lambda}}$ denotes the canonical projection of
$\tilde{\Lambda}_0\times V_0$ on $\tilde{\Lambda}$. By Proposition
\ref{r2.5}, $v(A,s)$ is a p.s.h. function or is equal $-\infty$
identically. But $v(0,s_0)=\langle
\mathcal{R},\pi_{\tilde{\Lambda}},(0,s_0)\rangle(\tilde{u})=\langle
R, \pi ,s_0\rangle(u)\ne -\infty$ then $v$ is a p.s.h. function.

Define $\pi_A:\Lambda_0\times V_0 \to V$ by $\pi_A(s,z):=s+Az$. Then
$\pi_A$ is a linear projection of $\Lambda_0\times V_0$ on
$\Lambda$. For $s\in \Lambda_0$, we have
$$\langle R, \pi_{A},s\rangle(u)=v(A,s).$$
Then $s\mapsto \langle R, \pi_{A},s\rangle(u)$ is a p.s.h. function
of $s\in\Lambda_0$ for all $A \in B(0,r)\backslash \mathcal{E}$,
where $\mathcal{E}\subset B(0,r)$ is a pluripolar set. This implies
that $s\mapsto\langle R, \pi_{A},s\rangle(u)$ is locally integrable
in $\Lambda_0$. Let $K$ be a compact subset of $\Lambda_0\times
V_0$. Then by the formula (\ref{4}), for $A\in
B(0,r)\backslash\mathcal{E}$, we have
\begin{eqnarray}\label{6}
\norm{uR\wedge \pi_{A}^* (i\d s_1\wedge
\d\bar{s}_1\wedge\ldots\wedge i\d s_m \wedge \d \bar{s}_m)}_K<
\infty.
\end{eqnarray}

We can obtain a strictly positive form on $K$ by taking a
combination of the forms $\pi_{A}^*(i\d s_1\wedge
\d\bar{s}_1\wedge\ldots\wedge i\d s_m \wedge \d \bar{s}_m)$ for
$A\in B(0,r)\backslash \mathcal{E}$. By the inequality (\ref{6}),
hence $uR$ has locally finite mass in $\Lambda_0\times V_0$. This
implies the theorem. We define $\ddc u\wedge R:=\ddc (uR)$, see e.g
[FS2].

\end{proof}
\end{appendix}

\textit{Acknowledgements}. I wish to thank T.C Dinh and N. Sibony
for their precious help during the preparation of this article.

\noindent Ngoc-Mai Pham

\noindent Math\'{e}matique- B\^{a}timent 425, UFR 8628

\noindent Universit\'{e} Paris-Sud, 91405 Orsay, France

\noindent {\tt Ngocmai.Pham@math.u-psud.fr}

\end{document}